\documentclass{article}

\usepackage{latexsym,url,epsf,epsfig}
\usepackage{amssymb,amsmath,graphicx,enumerate}

\newcommand{\cM}{\mathcal{M}}
\newcommand{\cC}{\mathcal{C}}
\newcommand{\cI}{\mathcal{I}}

\makeatletter
\def\GrabProofArgument[#1]{ (#1): \egroup\ignorespaces}
\def\proof{\noindent\textbf\bgroup Proof%
           \@ifnextchar[{\GrabProofArgument}{: \egroup\ignorespaces}}

\makeatother

\newcommand{\qed}{\hfill$\square$}

\newcommand{\N}{{\mathbb N}}

\let\emptyset=\varnothing
\newcommand{\old}[1]{{}}

\newtheorem {Theorem}{Theorem}
\newtheorem {Lemma}[Theorem]{Lemma}

\newtheorem {Proposition}[Theorem]{Proposition}
\newtheorem {Corollary}[Theorem]{Corollary}

\begin{document}

\title{Characterizing Matchings as \\
       the Intersection of Matroids\thanks{
An earlier version appears as an extended abstract
in the Proceedings of COMB'01~\cite{FFS01}.}}

\author{S\'{a}ndor P.\ Fekete\thanks{
Abteilung f\"ur Mathematische Optimierung, TU Braunschweig,
{\tt s.fekete@tu-bs.de}}
\and 
Robert T. Firla\thanks{
Institut f\"ur Mathematische Optimierung,
Otto-von-Guericke-Universit\"at Magde\-burg,
{\tt [firla,spille]@imo.math.uni-magdeburg.de}.
Supported by the ``Gerhard-Hess-For\-schungs\-f{\"o}rderpreis'' (WE
1462) of the German Science Foundation (DFG) awarded to R.~Weismantel.}
\and
Bianca Spille{\footnotemark[2]}
}

\date{}
\maketitle

\begin{abstract}
This paper deals with the problem of representing the matching
independence system in a graph as the intersection 
of finitely many matroids.
After characterizing the graphs for which the matching
independence system is the intersection of two matroids,
we study the function $\mu(G)$, which is the minimum
number of matroids that need to be intersected in order to
obtain the set of matchings on a graph $G$, and examine
the maximal value, $\mu(n)$, for graphs with $n$ vertices.
We describe an integer programming formulation for
deciding whether $\mu(G)\leq k$.
Using combinatorial arguments,
we prove that $\mu(n)\in \Omega(\log \log n)$. 
On the other hand, we establish that $\mu(n)\in O(\log n/\log \log n)$.
Finally, we prove that $\mu(n)=4$ for $n=5,\ldots,12$, and
sketch a proof of $\mu(n)$=5 for $n=13,14,15$.
\end{abstract}

{\small

}

\section{Introduction}\label{sec:intro}

Many combinatorial optimization problems can be viewed as optimization
problems over independence systems.
Some of them are polynomially solvable, e.g., spanning trees in
graphs, the branching problem in digraphs, or the matching problem.
Others are known to be NP-complete, like the traveling 
salesman problem (TSP) or the stable set problem 
(cf.~\cite{CCPS,GGL95,GLS88}).
Among the problems with a polynomial-time algorithm, the matching 
problem is generally recognized as one of the ``hardest'', 
and the famous blossom algorithm by Edmonds~\cite{Ed65} 
is one of the highlights of combinatorial optimization.

Another seminal result on the optimization in independence
systems was also given by Edmonds
\cite{Ed70}, who proved that the optimization problem 
over the intersection of two matroids is solvable in polynomial time. 
Algorithms for this problem were given by Edmonds \cite{Ed79}, Frank
\cite{Fr81}, and Lawler \cite{La75,La76}.
Unfortunately, this cannot be generalized to the case of three or more
matroids: As the NP-complete TSP can be written as
an optimization problem over the intersection of three matroids, 
it is highly unlikely that a polynomial-time algorithm exists.

In this paper, we consider a problem that combines aspects
of both problems that were solved by Edmonds:
How many matroids need to be intersected to characterize the matchings
of a graph? This problem is somewhat related to work by
Jenkyns~\cite{Jen} and Korte and Hausmann~\cite{KH78},
who described approximation guarantees for the greedy algorithm;
these guarantees depend on
the rank quotient, which depends on the number of matroids needed
for characterizing the independence system.
Clearly, we are not primarily interested in these approximation guarantees
for matching. However, we believe that the problem
of describing an independence system as the intersection of few matroids
is an interesting combinatorial problem in its own right.

The rest of this paper is organized as follows.
After some technical preliminaries in Section~\ref{preliminaries},
Section~\ref{main} gives a number of general structural results.
In particular, we give a precise characterizations of graphs
for which the set of matchings can be represented as the intersection
of two matroids, and an Integer Programming formulation
for the problem of minimizing the number of matroids that
are necessary for representing the matchings of an input graph.
In Section~\ref{sec:lower} we prove that $K_n$ needs at
least $\Omega(\log\log n)$ matroids, while Section~\ref{sec:upper}
establishes an upper bound of $O(\log n/\log\log n)$
for the number of matroids needed for any graph with $n$
vertices. The final Section~\ref{sec:tight} describes
the actual values of matroids that are necessary
for graphs with up to 15 vertices.

\section{Preliminaries}\label{preliminaries}

Let $S$ be a finite set and ${\cI}$ be a family of subsets of $S$.
${\cI}$ is an {\it independence system} on $S$ 
if $\emptyset\in{\cI}$ and 
if $J'\subseteq J$ and $J\in{\cI}$ then $J'\in{\cI}$.
The subsets of $S$ belonging to ${\cI}$ are called {\it independent},
otherwise {\it dependent}.
The minimal dependent subsets of $S$ are the {\it circuits} of ${\cI}$.
The {\it circuit system} ${\cC}$ of ${\cI}$ is the set of circuits of
${\cI}$ and 
${\cI}=\{J\subseteq S: C\not\subseteq J\mbox{ for all }C\in{\cC}\}$.
A maximal independent subset of a set $A\subseteq S$
is a {\it basis} of $A$.
An independence system ${\cI}$ on $S$ is a {\it matroid} 
if for every subset $A\subseteq S$ all its bases have the same
cardinality.
For further background, see Oxley \cite{Ox92} and Welsh \cite{We95}. 
Here we just state another useful fact.

\begin{Proposition}\label{circuitsystemmatroid}
  The circuit system ${\cal C}$ of an independence system ${\cal I}$
  is the circuit system of a matroid if and only if
  for all $C_1\neq C_2\in{\cal C}$ with $C_1\cap C_2\neq\emptyset$
  and for all $c\in C_1\cap C_2$
  there exists $C\in{\cal C}$ such that
  $C\subseteq C_1\cup C_2\setminus\{c\}$.
\end{Proposition}

Any independence system is the intersection of finitely many matroids:
Let $\cC$ be the set of circuits of $\cI$. For $C\in\cC$,
let $\cM_C$ be the matroid on $S$ with circuit system $\{C\}$,
i.e., $\cM_C=\{J\subseteq S: C\not\subseteq J\}$. Then 
$ {\cI}=\bigcap\{\cM_C: C\in{\cC}\}$.
This, however, may not be the most economical way to
describe ${\cI}$, because we may be able to
cover several circuits by the same matroid. 
In the following, we 
write $\mu({\cI})$ for the minimum
number of matroids necessary for this task.
Throughout the rest of this paper, ${\cI}$ is the set of 
matchings of a graph, which we describe in the following.

Consider a finite graph $G=(V,E)$. A {\it matching} in $G$ is a set
of edges that are pairwise disjoint. 
The set $M(G)$ of matchings in $G$
forms an independence system on $E$.
For simplicity, we write $\mu(G)$ for $\mu(M(G))$, i.e.,
\[
    \mu(G)=\min\{m\in\N: M(G) \mbox{ is the intersection of $m$
    matroids}\}.
\]
Furthermore, $\mu(n)$ is used for the maximum $\mu(G)$
on graphs with $|V|\leq n$, i.e.
\[ 
  \mu(n)=\max\{\mu(G): |V|\leq n\}.
\]

It is easy to see that the circuits of $M(G)$ are the sets 
that consist of two intersecting edges. (The reader should keep
in mind that throughout the rest of this paper, 
the term {\em circuit} refers to such a pair of edges.)
We call a circuit an {\it $i$-circuit} 
if its edges intersect in vertex $i$. 
We denote the circuit $\{ij,ik\}$ 
with the two  edges $ij$ and $ik$ by $i^{jk}$.

\section{General Characterizations}\label{main}
The following easy lemma implies that $\mu(n)=\mu(K_n)$.
\begin{Lemma}\label{subgraph}
  Let $G'=(V',E')$ be a subgraph of $G=(V,E)$.
  Then $\mu(G')\leq \mu(G)$.
\end{Lemma}

\proof
Let $M(G)$ be the intersection of matroids $\cM_1,\ldots,\cM_m$ on $E$.
For $a=1,\ldots,m$,
let $\cM_a':=\{J:J\subseteq E', J\in \cM_a\}$
be the restriction of $\cM_a$ to $E'$.
Then $M(G')$ is the intersection of the matroids
$\cM_1',\ldots,\cM_m'$ on $E'$.\qed\medskip

As a consequence of this lemma, the number of matroids needed to
represent the matchings in the complete graph $K_n$ on $n$ vertices is a
natural upper bound for the number of matroids needed to represent
the matching independence system of any graph on at most $n$ vertices,
i.e., $\mu(n)=\mu(K_n)$.

\subsection{Matchings as the Intersection of Two Matroids}\label{2matroids}

We present a complete characterization of the graphs 
for which the set of matchings is the
intersection of at most two matroids, by generalizing the concept of
bipartite graphs. 

For a bipartite graph $G=(V_1\times V_2,E)$, 
the set of matchings $M(G)$ is the intersection of 
two (partition) matroids.
More generally, we get the following.
\begin{Theorem}\label{th:m-partite}
 Let $G=(V,E)$ be an $m$-partite graph. 
 Then $M(G)$ is the intersection of at most
 $m$ matroids on $E$.
\end{Theorem}
\proof
Let $V_1\times\ldots\times V_m$ be an $m$-partition of $G$.
Use the partition matroids
\[
  \cM_i:=\{J\subseteq E : |J\cap \delta(v)|\leq 1 
                   \mbox{ for all } v\in V_i\}, \quad i=1,\ldots,m.\ \ \ \square
\]
\medskip

As  we will see in Section~\ref{sec:upper},
this upper bound has quite a bit of slack for large $\mu(G)$.
Moreover, there are non-bipartite graphs $G$ with $\mu(G)=2$,
as can be seen from the following characterization.
\begin{Theorem}\label{th:2mat}
  The set of matchings $M(G)$ of a graph $G=(V,E)$ 
  is the intersection of two matroids
  if and only if $G$ contains no odd cycle of cardinality $\geq 5$
  and each triangle of $G$ has at most one vertex with degree $> 2$.
\end{Theorem}
\proof
  {\bf a)}
  Let $G$ be a graph that contains no odd cycle of cardinality $\geq 5$
  and all triangles of $G$ have at most one vertex with degree $> 2$.
  We call a triangle {\it isolated} if all its vertices have degree $2$.
  Let $G'$ be the graph that we obtain from $G$ 
  by contracting all isolated triangles and 
  by deleting from any other triangle the edge  
  that connects the two vertices of degree $2$. 
  Then $G'$ is a bipartite graph with bipartition $V_1'\times V_2'$.
  For $a=1,2$, let $V_a$ consist of all vertices $v\in V$ that correspond
  to vertices of $V_a'$ in $G'$.
  Let ${\cal C}_a$ consist of all circuits $i^{jk}$ of $M(G)$
  such that $i\in V_a$ or $\{i,j,k\}$ is a triangle of~$G$ and 
  $j,k\in V_a$.    
  Then ${\cal C}_a$ is the circuit system of the matroid
  $
      \cM_a:=\{J\subseteq E\mid 
             C\not\subseteq J \mbox{ for all }C\in{\cal C}_a\}
  $
  and
  $M(G)$ is the intersection of $\cM_1$ and $\cM_2$.

  {\bf b)}
  Let $M(G)$ be the intersection of two matroids $\cM_1$ and
  $\cM_2$ on $E$.
  Suppose $G$ contains an odd cycle $v_1,v_2,\ldots,v_{2k+1}$, $k\geq 2$.
  Then 
  $\{v_1v_2,v_2v_3\}$, $\{v_2v_3,v_3v_4\}$, $\ldots,$ 
  $\{v_{2k+1}v_1,v_1v_2\}$ are circuits of $M(G)$ and hence circuits in
  at least one of the matroids. 
  W.l.o.g.~we obtain that $\{v_1v_2,v_2v_3\}$ and $\{v_2v_3,v_3v_4\}$
  are circuits of $\cM_1$, in contradiction to $\{v_1v_2,v_3v_4\}$
  being a matching, see Proposition~\ref{circuitsystemmatroid}.
  Suppose $G$ contains a triangle $\{u,v,w\}$ with two additional edges
  $uz$ and $vz'$ (possibly $z=z'$). 
  There are three $u$-circuits $\{uv,uw\}$, $\{uv,uz\}$,
  and $\{uw,uz\}$. W.l.o.g.~all three are circuits of $\cM_1$. Because
  $\{vw,uz\}$ is a matching it follows that $\{uw,vw\}$,
  $\{uv,vw\}$ are circuits of $\cM_2$ and hence also $\{uv,uw\}$. 
  But because $\{uw,vz'\}$ is a matching, $\{uv,vz'\}$ cannot be a
  circuit in any of the two matroids, a contradiction.  
  Consequently, $G$ contains no odd cycle of cardinality $\geq 5$ and
  each triangle of~$G$ has at most one vertex with degree $>2$.
  \qed

\subsection{IP-Formulation}\label{subsec:IP}
Next we describe a characterization of the problem whether the set
of matchings in a graph $G=(V,E)$ on $n$ vertices
can be represented as the intersection
of at most $m$ matroids in terms of necessary and sufficient conditions.
This characterization leads in a natural way to an IP-formulation of
the introduced problem, which can be solved by standard IP-solvers,
for at least  not too large values of $n$ and $m$.

Suppose first that the set of matchings $M(G)$ of $G$ is the
intersection of $m$ matroids $\cM_1,\cM_2,\ldots,\cM_m$ on $E$.
Any matching of $G$ must be independent in each of these matroids
and any circuit of $M(G)$ must be dependent (and hence a circuit)
in at least one of these matroids.
For any matroid $\cM_a$ and any circuit $\{ij,ik\}$ of $M(G)$ with $j<k$,
we introduce a $0/1$-variable $x^a_{ij,ik}$ 
which is $1$ if the circuit $\{ij,ik\}$ is dependent 
(and hence a circuit) in $\cM_a$ or $0$ otherwise, i.e.,
\[
    x^a_{ij,ik}
    =\left\{ \begin{array}{r@{\quad:\quad}l}
                    1 & \{ij,ik\}
                        \mbox{ is dependent in }\cM_a, \\
                    0 & \mbox{otherwise.}
             \end{array}
     \right.
\]

{\bf Cover condition.}
Because any circuit is dependent in at least one of the matroids,
we obtain the following {\it cover-inequalities}
\begin{equation}
    \sum_{a=1}^m x^a_{ij,ik} \geq 1
    \quad\mbox{for all } ij, ik\in E, j< k.\label{coverineq}
\end{equation}
  
{\bf Claw condition.}
For any $i,j,k,l$ different, it is not possible that exactly two of
the three circuits $\{ij,ik\}$, $\{ij,il\}$, and $\{ik,il\}$ of $M(G)$
are circuits in the same matroid $\cM_a$, i.e., we have
\[
       x^a_{ij,ik} +  x^a_{ij,il} + x^a_{ik,il} \not= 2
       \quad\mbox{for all } ij, ik, il\in E,j<k<l,\mbox{ for all }a .
\] 
This is modeled by the following {\it claw-inequalities}:
\begin{equation}
\renewcommand{\arraystretch}{1.1}
  \begin{array}{*6{@{\:}r}clll}\label{claw}
    + & x^a_{ij,ik} & + & x^a_{ij,il} & - & x^a_{ik,il} & \leq &  1 
    && \mbox{for all }ij, ik, il\in E,\\
    + & x^a_{ij,ik} & - & x^a_{ij,il} & + & x^a_{ik,il} & \leq &  1 &
 \smash{\Biggr\}} & j<k<l, \\
    - & x^a_{ij,ik} & + & x^a_{ij,il} & + & x^a_{ik,il} & \leq &  1 
      && \mbox{for all }a
  \end{array}
\end{equation}

{\bf Triangle condition.}
For any $ij, ik, jk\in E$ different, it is impossible that exactly two of
the three circuits $\{ij,ik\}$, $\{ji,jk\}$, and $\{ki,kj\}$ of $M(G)$
are circuits in the same matroid $\cM_a$, i.e., we have
\[
           x^a_{ij,ik} +  x^a_{ji,jk} + x^a_{ki,kj} \not= 2
           \quad\mbox{for all }ij, ik, jk\in E,i<j<k,\mbox{ for all }a.
\]
We obtain the {\it triangle-inequalities}:
\begin{equation}\label{triangle}
\renewcommand{\arraystretch}{1.1}
  \begin{array}{*6{@{\:}r}clll}
    +& x^a_{ij,ik} & + &  x^a_{ji,jk} & - & x^a_{ki,kj} & \leq & 1
    && \mbox{for all }ij, ik, jk\in E,  \\
    +& x^a_{ij,ik} & - &  x^a_{ji,jk} & + & x^a_{ki,kj} & \leq & 1 &
    \smash{\Biggr\}} & i<j<k, \\
    -& x^a_{ij,ik} & + &  x^a_{ji,jk} & + & x^a_{ki,kj} & \leq & 1
    && \mbox{for all }a
  \end{array}
\end{equation}

{\bf Matching condition.}
For any $i,j,k,l$ different with $ij,kl\in E$, 
$\{ij,kl\}$ is a matching in $G$. Hence,
it is not possible that both
circuits $\{ij,ik\}$ and $\{ki,kl\}$ of $M(G)$ are circuits 
in the same matroid $\cM_a$. 
This leads us to the {\it matching-inequalities}:
\begin{equation}\label{matching}
\renewcommand{\arraystretch}{1.1}
  \begin{array}{*3{@{\:}r}cl@{\quad}l}
    x^a_{ij,ik} & + & x^a_{ki,kl} & \leq & 1 
    & \mbox{for all }  ij,ik,kl\in E,j<k,i<l,j\neq l,\mbox{ for all }a\\
    x^a_{ik,ij} & + & x^a_{ki,kl} & \leq & 1 
    & \mbox{for all } ij,ik,kl\in E,k<j,i<l,j\neq l,\mbox{ for all }a\\
    x^a_{ik,ij} & + & x^a_{kl,ki} & \leq & 1 
    & \mbox{for all } ij,ik,kl\in E,k<j,l<i,j\neq l,\mbox{ for all }a
  \end{array}
\end{equation}

\begin{Theorem}\label{IPsolver}
 The set of matchings $M(G)$ of a graph $G$ 
 is the intersection of at most $m$ matroids 
 if and only if there exists $x^a_{ij,ik}\in\{0,1\}$
 for all $ij,ik\in E,j<k$, for all $a=1,\ldots,m$ such that all the
 inequalities in (\ref{coverineq}), (\ref{claw}), (\ref{triangle}),
 (\ref{matching}) are satisfied.
\end{Theorem}
 \proof
 We have to show that any feasible solution of the given system 
 leads to $m$ matroids such that $M(G)$ is their intersection. 
 Let $(x^a_{ij,ik})_{i,j,k,a}$ be such a feasible solution.
 For $a=1,2,\ldots,m$, define
 \[
   {\cal C}_a:=\{\{ij,ik\}: ij,ik\in E, j<k, x^a_{ij,ik}=1\}.
 \] 
 Due to the claw-, triangle-, and matching-inequalities, 
 ${\cal C}_a$ is the circuit system of a matroid.
 Its associated matroid is
 $
 \cM_a=\{J\subseteq E: C\not\subseteq J 
       \mbox{ for all } C\in{\cal C}_a\}.
 $
  We claim that $M(G)$ is the intersection of $\cM_1,\cM_2,\ldots,\cM_m$.

 Let $J\in M(G)$.
 Because any circuit of $\cM_a$ is a circuit of $M(G)$,  
 $J$ is an element of $\cM_1\cap\cdots\cap\cM_m$.
 Now let $J$ be an element of this intersection.
 Suppose $J$ is not in $M(G)$. 
 Then there is some circuit $\{ij,ik\}$ $(j<k)$ in $M(G)$ which is contained
 in $J$. 
 Due to the cover-inequalities, $\{ij,ik\}$ is a circuit in at least
 one matroid $\cM_a$, 
 in contradiction to $J$ being $\cM_a$-independent.
 \qed\medskip

For most IP-solvers it is
more efficient to solve an optimization problem 
instead of solving a feasibility problem.
We transform the feasibility problem 
into an optimization
problem by introducing additional $0/1$-variables $y_{ij,ik}$ 
for any circuit $\{ij,ik\}$ $(j<k$) of $M(G)$.
We replace the cover-inequalities (\ref{coverineq}) by the inequalities
\[
       \sum_{a=1}^m x^a_{ij,ik} \,-\, y_{ij,ik} \geq 0
       \quad\mbox{for all } ij,ik\in E, j< k
\] 
and try to maximize the sum of the new  $y$-variables, i.e.,
$ \max\sum_{ij,ik\in E, j< k}y_{ij,ik}.$
This means we want to cover as many circuits as possible.
Consequently, the original feasibility problem has a feasible solution  
if and only if the new program has a feasible solution
in which all $y$-variables are equal to $1$.
Note that the $0$-vector is a feasible starting solution for this
integer program.

Nevertheless, these problems are still hard to solve.
The integer programs that we explore are quite large
and grow very fast because the problems have  
${\cal O}(mn^3)$ variables and  ${\cal O}(mn^4)$ constraints.
Therefore, current IP-solvers, e.g., 
CPLEX\footnote{CPLEX Linear Optimizer 6.0 (with Mixed Integer \&
  Barrier Solvers); \copyright\, ILOG Inc., Incline Village, NV, USA.} 
or SIP\footnote{Solving Integer Programs 1.1 by Alexander Martin (TU
  Darmstadt, Germany), unpublished.}, 
are unable to handle them in reasonable time even for moderate 
values of $n$ and $m$ (e.g., $G=K_n$ with $n=13$, $m=4$).

\section{Lower Bounds}\label{sec:lower}

In the following, we use the notation $\nu(m)$ to indicate the largest
$n$ for which $\mu(n)\leq m$, i.e.,
\[
    \nu(m)=\sup\{n\in\N:\mu(n)\leq m\}.
\]
The following result shows that 
$\nu(m)$ is indeed finite and grows at most doubly exponentially.

\begin{Theorem}\label{th:lower}
  $\nu(m) \leq 2^{2^{3m}-1}-1$.
\end{Theorem}

\proof
Let $G=(V,E)$ be the complete graph $K_n$ on $n$ vertices
and $\mu(n)\leq m$.
We start by introducing some technical terms.
For a vertex $i$, consider all directed edges $ij, j\in V\setminus\{i\}$.
(The difference between $ij$ and $ji$ is only important for the 
following definition.)
If an $i$-circuit involving edge $ij$ is covered by matroid
$\cM_a$, we say that $ij$ has {\em color} $a$.
The set of colors of an edge defines its {\em color type}.
The {\em color class} of any vertex $i$ is the
set of color types of all the edges $ij$, $j\in V\setminus\{i\}$.
Finally, a {\em degenerate triangle} consists of three different vertices
$i,j,k$, such that all three circuits $\{ij,ik\}$,
$\{ji,jk\}$, $\{ki,kj\}$ are covered by the same matroid, say, $\cM_a$.
For such a triangle, it is an easy consequence of the matching condition
described in Section~\ref{subsec:IP} that none of the
other circuits involving the edges $ij,ik,jk,ji,ki,kj$ may be covered
by $\cM_a$; i.e., for $l\not\in\{i,j,k\}$,
none of the circuits 
$\{ij,il\}$, $\{ik,il\}$, $\{ji,jl\}$, $\{jk,jl\}$,
$\{ki,kl\}$, $\{kj,kl\}$ may be covered by $\cM_a$ as well.

We start by eliminating covers with degenerate triangles.
If there are any degenerate triangles in a circuit cover
with $m$ matroids, we can construct a circuit
cover with $3m$ matroids that is free of degenerate triangles: 
For each degenerate triangle covered by
matroid $\cM_a$, cover one circuit by matroid $\cM_a$,
the other two by the additional matroids $\cM_{a'}$ and $\cM_{a''}$.
By the observation at the end of the preceding paragraph,
this does not affect any other circuits already covered by $\cM_a$, 
so all conditions described in Section~\ref{subsec:IP} are still valid.
Therefore, this yields indeed a feasible set of matroids.

Now consider the situation in the absence of degenerate triangles.
As any (directed) edge $ij\in E$ is part of some circuit, and each circuit is
covered by some matroid, there can be at most $2^m-1$ different color types
for $m$ different matroids.
Moreover, there are at most $2^{2^m-1}-1$
different color classes of vertices.
Furthermore, any $i$-circuit formed
by a pair of edges $ij$, $ik$ with $j\neq k$
must be covered by some matroid, so no two color types
in any valid color class can be disjoint. 

Assume that we have two vertices ($i$ and $j$) of the same 
color class. Then their connecting edge $ij$ is part of 
a circuit $\{ij,ik\}$ that is covered by some matroid $\cM_a$;
at the same time, there must be a circuit $\{ji,jl\}$
that is also covered by $\cM_a$. If $k\neq l$, we get a violated
matching condition. If $k=l$, we use the fact that $\{i,j,k\}$
is not a degenerate triangle, and we conclude that we get a violated
triangle condition.

Therefore, there can be at most $2^{2^{m}-1}-1$ vertices in the absence
of degenerate triangles, and not more than $2^{2^{3m}-1}-1$ in general.
\qed\medskip

Quite clearly, upper bounds for $\nu(m)$ correspond to lower bounds
for $\mu(n)$. In particular, we get
\begin{Corollary}\label{cor:lower}
  $\mu(n)\in \Omega(\log \log n)$.
\end{Corollary}

\section{Upper Bounds}\label{sec:upper}

Using a recursive construction, we can
show that $m$ matroids suffice to generate the matchings
of any graph with $O(\sqrt[3]{m!})=2^{O(m\log m)}$ vertices.
As above, this yields an upper bound for $\mu(n)$.
More precisely, we show
\begin{Theorem}\label{uppermatching}
\[\nu(m)\geq 3\prod_{i=0}^{\lceil\frac{m}{3}\rceil-2} (m-3i).\]
\end{Theorem}
\proof
We proceed by induction. Clearly, the claim is true for $m=1,2,3$, as 
all three circuits of $K_3$ can be covered by the circuit system
of just one matroid.
Now suppose the claim was true for $m-3$.
Consider $m$ vertex sets $V_1,\ldots,V_m$, each consisting of 
$\nu(m-3)$ vertices.

For proving the overall claim, it suffices
to describe circuit systems $\cC_1,\ldots,\cC_m$
that satisfy the following conditions. 
\begin{enumerate}
 \itemsep-.5ex
 \item[(1)] All circuits $v^{uw}$ with $v,u,w\in V_k$ are contained 
   in some $\cC_i$.
 \item[(2)] All circuits $v^{uw}$ with $v\in V_k$, $u,w\not\in V_k$ 
   are contained in some $\cC_i$.
 \item[(3)] All circuits $v^{uw}$ with $v,u\in V_k$, $w\not\in V_k$ 
   are contained in some $\cC_i$. 
 \item[(4)] Each resulting $\cC_i$ is the circuit system of a matroid.
\end{enumerate}

\begin{figure}
\centerline{\epsfig{file=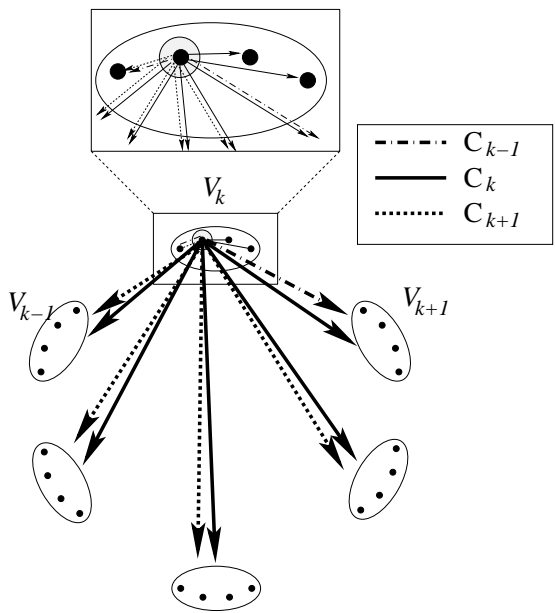,width=0.6\textwidth}}  
\caption{The structure of those circuits in $\cC_{k-1},\cC_k,\cC_{k+1}$
that have a central vertex in $V_k$.}
\label{fig:picture}
\end{figure}

In the following, indices are taken modulo~$m$.
By the induction hypothesis, the set of matchings on $V_k$ is the
intersection of $m-3$ matroids. Denote the corresponding circuit systems
by $\cC_1^k,\ldots,\cC_{m-3}^k$.
For $i=1,\ldots,m$ define
\[
     \cC_i':=\cC^{i+2}_1\cup \cC^{i+3}_2\cup\ldots\cup\cC^{i-2}_{m-3}
     =\bigcup_{j=1}^{m-3}\cC_j^{i+1+j}.
\]
In particular, $\cC_{k-i-1}^k\subseteq \cC_i'$ is the set of all 
circuits within $V_k$ that are covered by~$C_i'$.
Then $\cC_i'$ is a circuit system on $E$ and $\cC_1',\ldots,\cC_m'$
satisfy condition (1).
The circuit systems $\cC_1,\ldots,\cC_m$ arise now by adding
$\cC_i'$ to $\cC_i$ and for $k=1,\ldots,m$:
\begin{itemize}
  \item Add $\{v^{uw}:v\in V_k, u,w\not\in V_k\}$
        to $\cC_k$. \\
        Add $\{v^{uw}:v\in V_k, u,w\in V_{k+1}\}$
        to $\cC_{k-1}$.\\
        Add $\{v^{uw}:v\in V_k, u,w\not\in V_k\cup V_{k+1}\}$
        to $\cC_{k+1}$. 
  \item Now consider $V_k=\{v_{(k-1)\nu(m-3)+1},\ldots,v_{k\nu(m-3)}\}$. \\
        Add the following circuits to $\cC_k$
        \[
                 \{v_i^{v_j w}: i<j,v^i,v^j\in V_k,w\not\in V_k\}\cup
                 \{v_i^{v_j v_l}:i<j,l, v^i,v^j,v^l\in V_k\}.
        \]
        Add the following circuits to $\cC_{k-1}$
        \[
                 \{v_i^{v_j w}:i>j,v^i,v^j\in V_k,w\in V_{k+1}\}\cup
                 \{v_i^{v_j v_l}:i>j,l, v^i,v^j,v^l\in V_k\}.
        \]
        Add the following circuits to $\cC_{k+1}$
        \[
                 \{v_i^{v_j w}:i>j,v^i,v^j\in V_k,
                                   w\not\in V_k\cup V_{k+1}\}\cup
                 \{v_i^{v_j v_l}:i>j,l, v^i,v^j,v^l\in V_k\}.
        \]
\end{itemize}
See Fig.~\ref{fig:picture} for an illustration.
There, arrowheads point away from the central vertices of circuits.

The resulting circuit systems $\cC_1,\ldots,\cC_m$
satisfy conditions (1), (2), and (3). 
In order to see that each $\cC_i$ is the circuit system of a matroid,
we need to verify the claw-, triangle-, and matching-conditions 
described in Section~\ref{subsec:IP}.

{\bf Claw condition:}
We need to verify that for any 
two circuits $v^{uw_1}$, $v^{uw_2}$ that are both contained in 
some $\cC_i$, $v^{w_1w_2}$ is also contained in $\cC_i$. 

Let $v\in V_k$. 
If $i\not\in\{k-1,k,k+1\}$, then the involved circuits
must be completely within $V_k$
and they are all covered by $\cC^k_{k-i-1}\subseteq \cC_i$.
Otherwise, $i\in\{k-1,k,k+1\}$ and the property is easily verified.

{\bf Triangle and matching conditions:}
Let $v^{uw_1}\in\cC_i$ for some $i$.
We need to verify that if $u^{vw_2}\in\cC_i$ for some $w_2$
then $w_2=w_1$ and $w_1^{uv}\in\cC_i$.
  
Let $v\in V_k$. 
If $i\not\in\{k-1,k,k+1\}$, then the involved circuits
must be completely within $V_k$ and hence, being in $\cC_i$ means
being in $\cC^k_{k-i-1}$.
Because $\cC^k_{k-i-1}$ is the circuit system of a matroid, 
the triangle and matching conditions are satisfied.
In the following $i\in\{k-1,k,k+1\}$.
It suffices to show that 
no $u^{vw_2}$ is contained in $\cC_i$. 
If $u\not\in V_k\cup V_{k+1}$ then $i=k$ or $i=k+1$.
On the other hand, none of the circuits $u^{vw_2}$ is in $\cC_k$ or
$\cC_{k+1}$. 
If $u\in V_{k+1}$ it follows by our construction that $i=k$ or $i=k-1$ 
while $u^{vw_2}$ is in $\cC_{k+2}$ or $\cC_{k+1}$. 
Finally, we consider $u\in V_{k}$. 
Let $v=v_h$ and $u=v_j$.
If $h<j$ then $i=k$ while $u^{vw_2}$ is in $\cC_{k-1}$ or $\cC_{k+1}$.
If $h>j$ then $i=k-1$ or $i=k+1$ while $u^{vw_2}$ is in $\cC_k$.

\medskip
Consequently, $\cC_1,\ldots,\cC_m$ are circuit systems of matroids, 
say $\cM_1,\ldots,\cM_m$,
and all circuits of $M(G)$ are covered by some $\cC_i$.
Hence, $M(G)$ is the intersection of the matroids
$\cM_1,\ldots,\cM_m$.
Therefore, $\nu(m)\geq m\cdot\nu(m-3)$.
With $\mu(j)\geq 3$ for $j\in\{1,2,3\}$ we obtain
\[\nu(m)\geq 3\prod_{i=0}^{\lceil\frac{m}{3}\rceil-2} (m-3i).\]
\qed\medskip

Now the lower bound on $\nu(m)$ 
implies an upper bound on $\mu(n)$:

\begin{Corollary}\label{cor:upper}
  $\mu(n)\in O(\log n/\log \log n)$.
\end{Corollary}

\proof
For $n \geq s = 2^{c m \log m}$, i.e., $m \log m = \frac{\log s}{c}$,
we get $\log m + \log\log m=\log (m \log m) = \log \log s - \log c$.
This implies $2 \log m \geq \log \log s$ for sufficiently large $m$.
Therefore, 
$\frac{2\log n}{c\log\log n}\geq \frac{2\log s}{c\log\log s}\geq m$.
\qed

\section{Tight Bounds for $\mu(n)$}\label{sec:tight}

Lawler mentioned in \cite{La76} that the nonbipartite matching
problem can be formulated as an intersection problem involving two
partition matroids, but with additional constraints in the form of
symmetry conditions. Nevertheless, we can give an elementary proof
of the following.
\begin{Theorem}\label{th:mu5}
  $\mu(n)=4$ for $n=5,\ldots,12$.
\end{Theorem}
\proof
{From} Theorem~\ref{uppermatching} we obtain $\nu(4)\geq 12$,
hence, $\mu(n)\leq 4$ for $n\leq 12$.

We complete the proof by showing that $\mu(5)\geq 4$.
Let $G$ be the complete graph $K_5$ 
with vertex set $\{1,2,3,4,5\}$ and 
$M(G)$ the set of matchings of $G$.
Suppose $M(G)$ is the intersection of three matroids $\cM_1,\cM_2$,
and $\cM_3$ on $E$.

Consider the circuits that correspond to one vertex
with incident edges $a,b,c,d$. 
These are $\{a,b\}$, $\{a,c\}$, $\{a,d\}$, $\{b,c\}$, $\{b,d\}$, $\{c,d\}$.

Suppose exactly three of them are $\cM_1$-circuits,
w.l.o.g., $\{a,b\},\{a,c\},\{b,c\}$.
Without loss of generality, $\{a,d\}$ and $\{b,d\}$ 
and hence $\{a,b\}$ are $\cM_2$-circuits.
There exist edges $e$ and $f$ such that $\{a,b,e\}$ and $\{b,d,f\}$
are triangles:

\noindent
Because $\{c,e\}$ is a matching $\{a,e\}$ is no $\cM_1$-circuit and
because $\{d,e\}$ is a matching $\{a,e\}$ is no $\cM_2$-circuit,
hence, $\{a,e\}$ is an $\cM_3$-circuit,
Similarly we obtain $\{b,e\}$ is an $\cM_3$-circuit. 
Consequently, $\{a,b\}$ is an $\cM_1$-, an \mbox{$\cM_2$-,} 
and an $\cM_3$-circuit.
Because $\{a,f\}$ is a matching, 
there is no possibility for $\{b,f\}$ 
to be a circuit in one of the three matroids,
a contradiction.

Hence, it is not possible that any matroid has exactly three of the
six circuits that correspond to one vertex as circuits.
If a matroid has at least four of the
six circuits that correspond to one vertex as circuits then all six
circuits are its circuits.
Therefore, all six circuits are circuits in the same matroid or 
each of the three matroids has exactly two of the six circuits as
circuits, these must be disjoint.
Consequently,
there are $i_1,i_2,i_3,i_4,i_5\in\{1,2,3\}$
such that $\{12,15\}$ and $\{13,14\}$ are $\cM_{i_1}$-circuits, 
          $\{12,23\}$ and $\{24,25\}$ are $\cM_{i_2}$-circuits,
          $\{23,34\}$ and $\{13,35\}$ are $\cM_{i_3}$-circuits,
          $\{34,45\}$ and $\{14,24\}$ are $\cM_{i_4}$-circuits, and
          $\{45,15\}$ and $\{25,35\}$ are $\cM_{i_5}$-circuits.
Up to symmetry, we obtain $i_1=i_2$, 
i.e., $\{12,15\}$, $\{13,14\}$, $\{12,23\}$, and $\{24,25\}$ 
are $\cM_{i_1}$-circuits, in contradiction to $\{15,23\}$ being a matching. 
\qed\medskip

We can also give a positive result concerning matching and the
intersection of three matroids: 
\begin{Theorem}\label{4partite}
 Let $G=(V,E)$ be a $4$-partite graph.
 Then the set of matchings~$M(G)$ of~$G$ is the intersection of at most
 three matroids on $E$.
\end{Theorem}
 \proof
 Let $V_1\times V_2\times V_3\times V_4$ be a $4$-partition of $G$.
 For $a=1,2,3$, let $\cC_a$ consist of all $V_a$-circuits and of all 
 $V_4$-circuits $i^{jk}$ with $j,k\not\in V_a$ ($i\in V_4$).
 Then $\cC_a$ is the circuit system of a matroid. Its associated
 matroid is
 $\cM_a=\{J\subseteq E:C\not\subseteq J\mbox{ for all }C\in\cC_a\}$.
 It is easy to see that $M(G)$ is the intersection of the matroids
 $\cM_1,\cM_2,\cM_3$ on $E$.
 \qed\medskip

This implies that
the set of matchings of any subgraph of $K_5$
is the intersection of at most three matroids, i.e.,
$K_5$ is the smallest graph for which the set of matchings is
not the intersection of three matroids.

Theorem~\ref{uppermatching} implies that $\nu(5)\geq 15$,
hence, $\mu(n)\leq 5$ for $n=13,14,15$.
Using refined versions of the techniques for the lower bound
described in Section~\ref{sec:lower}, we can show the following:
\begin{Theorem}\label{th:mu13}
  $\mu(n)=5$ for $n=13,14,15$.
\end{Theorem}

The proof proceeds by showing that for $m=4$, there are at most 12 
``basic'' color classes, i.e., color classes that cannot be simplified
by deleting some of the circuits from some of the matroids.
This implies the claim in the absence of degenerate triangles,
as there cannot be any two vertices from the same basic color class.
Furthermore, any pair of vertices from the same basic color class
(which forces a degenerate triangle) eliminates another basic
color class; again, the claim follows.

As full details are rather tedious and probably not of sufficient
interest to the reader, they are omitted.

\end{document}